\documentclass[10pt]{amsart}
\usepackage{amssymb}
\usepackage{amsfonts}
\usepackage{amsmath}
\usepackage{amsthm}

\newtheorem{teo}{Theorem}
\newtheorem{prop}[teo]{Proposition}
\newtheorem{definicion}[teo]{Definition}

\newtheorem{lema}[teo]{Lemma}

\DeclareMathOperator{\codim}{codim}
\newcommand{\cuerpo}[1]{\mathbb{#1}}
\newcommand{\zP}{\cuerpo{P}}
\DeclareMathOperator{\rg}{rank}
\DeclareMathOperator{\se}{Sec}
\title{On the rank of a binary form}

\author{Gonzalo Comas} 
\address{Departamento de Matematica FCEyN-UBA
\\ Ciudad Universitaria (1428) Capital Federal
\\ Argentina}
\email{gcomas@dm.uba.ar}
\author{Malena Seiguer}
\email{mseiguer@dm.uba.ar}

\begin{document}

\begin{abstract}
We describe in the space of binary forms of degree $d$
the strata of forms having constant rank.
We also give a simple algorithm to determine the rank of a given form.
\end{abstract}

\maketitle

\begin{section}{Introduction}
Let $V$ be a vector space
of  dimension two over an
algebraically closed field $K$ of characteristic zero. Let $q$ be a binary 
form of degree $d$ over $V$.
We define the {\it rank} of 
$q$ (denoted $\rg q$) to be the least integer $r$ such that $q$ can be written as a sum 
of $r$ $d$th powers of linear forms:
$$
q=L_1^d +\dots + L_r^d.
$$
In this article we describe the sets of homogeneous forms of a  given rank.

We denote with $S_{d,r}$ the set of binary forms of degree $d$ and rank $r$.
We show that if a limit of forms in $S_{d,r}$ has rank greater than $r$, then
it will have rank
greater than or equal to $d-r+2$. More precisely, we prove that
$$
\overline{S_{d,r}} \setminus \overline{S_{d,r-1}}
=S_{d,r} \cup S_{d,d-r+2}
$$
for $2 \leq r \leq  \left[\frac d 2\right]+1$; and $\overline{S_{d,1}}=S_{d,1}$
is the set of $d$th powers of linear forms
(see theorem \ref{elteorema}).

We can restate the result of the theorem as follows
$$
\overline{S_{d,r}} = \overline{S_{d,r-1}} \cup S_{d,r} \cup S_{d,d-r+2}.
$$
Applying this result repeatedly, we get
a characterization of the possible ranks of degenerated forms
$$
\overline{S_{d,r}}=\bigcup_{i=1}^{r} S_{d,i} \cup \bigcup_{i=d-r+2}^{ d}
S_{d,i}.
$$

We also show how to determine the rank of a given form $q$
by computing the rank of an explicit matrix,
solving a linear system of equations,
and deciding if a polynomial constructed
from this solution has multiple roots.

{\bf Acknowledgments.} We wish to thank F. Cukierman for suggesting this problem
and for his encouragement. We also thank J.M. Landsberg for useful discussions.

\end{section}

\begin{section}{Rank of a binary form}

Let $V=K^2$ and let $S_1=V^\ast$ be the vector space of linear
polynomials in two variables $x,y$.
Also let $S_d$ be the space of homogeneous polynomials 
of degree $d$ in  variables $x,y$. 
We define the rank of a binary form $q \in S_d$ as follows
\medskip
\begin{definicion}
The rank of $q$
is the least integer $r$ such that $q$ can be written
as
$$
q=L_1^d +\dots + L_r^d
$$ 
where $L_i \in S_1$, $i=1,\dots,r$.
\end{definicion}

Let $X\subset \zP(S_d)$ be  the Veronese curve, that is, the image of the map
$$
\begin{array}{ccc}
\zP(S_1) & \rightarrow & \zP(S_d) \\
{[L]}    & \mapsto	 & {[L^d]}.
\end{array}
$$
Then the rank of a binary form $q$ is the least integer $r$ such that
$[q]$ lies on a linear space 
spanned by  $r$ elements lying on $X$. 
For instance, the forms of rank one are exactly those lying on $X$,
and those lying on a secant line to $X$ have rank less than or equal to 2.
However, not all forms lying on the secant variety to the curve have rank 
less than or equal
to 2: in fact, a form lying on a tangent line but not on the curve has rank $d$
(\cite{Ha}, pag. 147).
Therefore, this notion of rank does not behave upper-semicontinously.

Similarly, 
forms lying on a secant $k$-plane spanned by $k+1$ points
of the curve have rank less than or equal to $k+1$, but
their degenerations may have higher rank.

Let 
$$
S_{d,r}=\{q \in S_d: \rg q =r\}
$$
be the
set of forms of rank $r$.  
Let $\se^k(X)$ denote  the $(k+1)$-secant variety of $X$, that is,
the variety defined as the closure of the union of $k$-planes 
spanned by $k+1$ points in $X$.
Note that $\se^k(X)$ is the projectivization of  the closure of $S_{d,k+1}$,
that is, $\se^k(X)=\{[q] \in \zP(S_d) : q \in \overline{S_{d,k+1}}\}$.

We will prove the following theorem
\medskip
\begin{teo}\label{elteorema}
For each integer  $k$, $0 \leq k \leq \left[\frac d 2\right]$, we have
$$\overline {S_{d,k+1}} \setminus \overline {S_{d,k}}
=S_{d,k+1} \cup S_{d,d-k+1}
$$
where $\overline{S_{d,0}}=\emptyset$.

For $q \in \overline {S_{d,k+1}} \setminus \overline {S_{d,k}}$
we have $\rg q=k+1$ if and only if $[q]$
lies on a secant $k$-plane to $X$; otherwise, $\rg q=d-k+1$.
\end{teo}

\bigskip

For the proof, 
we will use the following alternative way of describing the Veronese curve.
Consider the map
$$
\begin{array}{ccc}
\zP^1(K) & \rightarrow & \zP(S_d^\ast) \\
{[\alpha]}    & \mapsto	 & {[ev(\alpha)]}
\end{array}
$$
where $ev(\alpha)$ is the linear functional 
given by evaluation of polynomials at $\alpha \in K^2$.
The Veronese curve is the image of this map. More precisely,
let $\{x,y\}$ be a basis for $S_1$, and $\{x^d,x^{d-1}y,\dots,xy^{d-1},y^d\}$
a basis for $S_d$. We consider in $S_d^\ast$
the base dual to that of $S_d$.
Then the image of a point in $\zP(S_1)$ with homogeneous
coordinates $[t,u]$  by the first map has homogeneous coordinates
$\left[ t^d , {\binom{d}{1}}t^{d-1}u, \dots,{\binom{d}{d-1}}tu^{d-1},u^d\right]$
in $\zP(S_d)$.
On the other hand, a point in $\zP^1(K)$ with homogeneous coordinates $[t,u]$
is mapped by the second map
to  a point in $\zP(S_d^\ast)$ with homogeneous coordinates
$\left[ t^d , t^{d-1}u, \dots,tu^{d-1},u^d\right]$.

Therefore, the isomorphism given by
$$
\begin{array}{ccc}
\zP(S_d^\ast) & \rightarrow & \zP(S_d) \\
{[Z_0,Z_1,\dots,Z_{d-1},Z_d]}    & \mapsto	 & {\left[Z_0,{\binom{d}{1}}Z_1,\dots,
{\binom{d}{d-1}}Z_{d-1},Z_d\right]}
\end{array}
$$
restricts to an isomorphism between the image of the second Veronese map
and the first one.

\medskip

Using this isomorphism, we define  the rank of a functional $\varphi \in S_d^\ast$
as follows
\medskip
\begin{definicion}
Let $\varphi \in S_d^\ast$ be a linear functional. The rank of $\varphi$
is the least integer $r$ such that $\varphi$ can be written as
$$
\varphi = ev(\alpha_1)+\dots+ev(\alpha_r)
$$ 
where $\alpha_i \in K^2$, $i=1,\dots,r$. 
\end{definicion} 
Notice that if $\varphi = \lambda_1 ev(\alpha_1)+\dots+ \lambda_r ev(\alpha_r)$,
then $\varphi$ is the sum of $ev(\alpha'_i)$, where $\alpha'_i=\omega_i\alpha_i$
and $\omega_i^d=\lambda_i$.  

\bigskip

We will need two results regarding secant varieties to the Veronese curve
(\cite{Ha}, propositions 11.32 and 9.7).

\medskip

\begin{prop}
Let $X \subset \zP(S_d^\ast)$ be the Veronese curve. The secant variety
$\se^k(X)$ has dimension $\min\{2k+1,d\}$.
\end{prop}

Note that if
$d=2n$ or $d=2n+1$, $\se^n(X)=\zP^d$.
This means that computing the rank of points lying on
$\se^k(X)$ for $k\leq\left[\frac d 2 \right]$ gives a complete stratification by rank.

\bigskip

Given  $[\varphi]\in \zP(S_d^*)$, we will
consider its homogeneous coordinates $[Z_0,Z_1,\ldots,Z_d]$
where $Z_i=\varphi(x^{d-i}y^{i})$.

\medskip

\begin{prop}\label{harris}
For any $k \leq s,d-s$, the rank $k$ determinantal variety associated to the 
matrix
$$
M=\left [
\begin{array}{cccc}
Z_0 & Z_1 & \dots & Z_{d-s}\\
Z_1 & Z_2 & \dots & Z_{d-s+1}\\
\vdots & \vdots &  & \vdots\\
Z_{s} & Z_{s+1} & \dots & Z_d\\
\end{array}
\right]$$
is the $k$-secant variety $\se^{k-1}(X)$ of the Veronese curve $X \subset \zP(S_d^*)$.

\end{prop} 

\bigskip

We will  also need Bertini's Theorem (\cite{Ha}, pag. 216)
\medskip
\begin{teo}[Bertini's Theorem]\label{bertini2}
If $X$ is any quasiprojective variety, $f:X \rightarrow \zP^n$
a regular map, $H \subseteq \zP^n$ a general hyperplane, and $Y=f^{-1}(H)$,
then 
$$
Y_{\text{sing}}=X_{\text{sing}}\cap Y.
$$
In particular, we will use the formulation: the general member of a linear
system on $X$ is smooth outside the singular locus of $X$
and the base locus of the linear system.
\end{teo}

\bigskip

We will need the following characterization of the secant variety $\se^k(X)$,
that will allow us to recognize which points lie on secant $k$-planes
and which are degenerate.

\medskip 		

\begin{lema}\label{secante}
Let $[\varphi] \in \zP(S_d^\ast)$ and let $k$ be an integer, 
$0 \leq k < \left[\frac d 2 \right]$. Then $[\varphi] \in \se^k(X)$
if and only if $\varphi(L_1L_2\dots L_{k+1} S_{d-k-1})=0$ 
for some linear forms $L_i$, $i=1,\dots ,k+1$.
\end{lema}

\begin{proof}
Consider the bilinear form $B:  S_{d-k-1}  \times S_{k+1} \rightarrow K$ given by 
$B(f,g)=\varphi(fg)$. The matrix of $B$ in basis 
$$ \{x^{d-k-1},x^{d-k-2}y,\dots,y^{d-k-1}\}
\hskip 10pt
\{x^{k+1},x^{k}y,\dots,y^{k+1}\} 
$$
is
$$
M=\left [
\begin{array}{cccc}
Z_0 & Z_1 & \dots & Z_{k+1}\\
Z_1 & Z_2 & \dots & Z_{k+2}\\
\vdots & \vdots &  & \vdots\\
Z_{d-k-1} & Z_{d-k} & \dots & Z_d\\
\end{array}
 \right]$$
where $Z_i=\varphi(x^{d-i}y^i)$ denote the coordinates of $[\varphi]$.
Using proposition \ref{harris},
 $[\varphi] \in \se^k(X)$ if and only if the rank of this matrix is 
less than or equal to $k+1$.
Equivalently, $[\varphi] \in \se^k(X)$
if and only if there exists a polynomial $g \in S_{k+1}$
such that $B(S_{d-k-1},g)=0$. Factoring $g$ as $g=L_1\dots L_{k+1}$
we get our result. 
\end{proof}

\bigskip

In the following lemma we give a necessary and sufficient condition
for the rank of  $\varphi$ to be
less than or equal to a given integer $r$.
We shall denote with $\Delta_r $
the discriminant hypersurface in the affine space $S_r$, defined as the 
locus of polynomials with multiple roots in $\zP^1(K)$. 
It is a well known fact that $\Delta_r$ is an irreducible
hypersurface.

\medskip

\begin{lema}\label{igual}
$\varphi=\sum_{i=1} ^r \lambda_i ev(\alpha_i)$, where  $[\alpha_i]$ are
distinct points in $\zP^1(K)$
if and only if the set
$A=\{f\in S_r : \varphi(fg)=0\;\forall g\in S_{d-r}\}$ 
is not contained in $\Delta_r$.

In particular, every  $\varphi$ is a linear combination
of $d$ elements in $X$, therefore $d$ is an upper bound for the rank of a binary form.
\end{lema}

\begin{proof}
Let $\varphi=\sum_{i=1} ^r \lambda_i ev(\alpha_i)$.
The polynomial $f$ with roots $[\alpha_1],\dots,[\alpha_r]$
lies on $A$, therefore $A$ is not contained in $\Delta_r$.

For the converse, let us consider a polynomial
$f=L_1.L_2\dots L_r \in A\setminus\Delta_r$, 
where $L_1,\dots,L_r$ are linear forms.
For each $i=1,\dots,r$, let $[\alpha_i]=[t_i,u_i] \in \zP^1(K)$ 
be the zero of $L_i$. Since $f \not \in \Delta_r$,
$[\alpha_1],\dots,[\alpha_r]$ are distinct.
We will show that
$\varphi=\sum_{i=1} ^r \lambda_i ev(\alpha_i)$
by proving that
$\{\varphi,ev(\alpha_1),\dots,ev(\alpha_r)\}$ 
is a linearly dependent set.
Consider the $(r+1)\times(d+1)$ matrix having as its rows  the coordinates
of these functionals in the base dual 
to $\{x^d,x^{d-1}y,\dots,xy^{d-1},y^d\}$.
$$
\left[
\begin{array}{cccc}
\varphi(x^d) & \varphi(x^{d-1}y) & \dots & \varphi(y^d)\\
t_1^d & t_1^{d-1}u_1 & \dots & u_1 ^d\\
\vdots & \vdots & \ddots & \vdots\\
t_r^d & t_r^{d-1}u_r & \dots & u_r ^d\\
\end{array}
\right]
$$
We claim that this matrix does not have maximal rank.
Let us consider the maximal minor  obtained by 
choosing $r+1$ columns in this matrix.
Using the linearity of 
$\varphi$, this minor can be expressed as
$\varphi(g)$, where $g\in S_d$ is a polynomial having 
$[\alpha_1],\dots,[\alpha_r]$ as roots in $\zP^1(K)$. Thus, $\varphi(g)=0$.

When $r=d$, $A$ is the kernel of $\varphi$, which is a hyperplane.
Since $\Delta_d$ is an irreducible hypersurface, $A \not \subseteq \Delta_d$,
and therefore  $\varphi$ is a linear combination 
of $d$ elements in $X$.
\end{proof}

In the next proposition we show a lower bound for
the increase of  rank when  passing to a degenerate position.  

\medskip

\begin{prop}\label{ida}
Let $k\leq\left[\frac {d-1} 2 \right]$.
If $[\varphi] \in \se^k(X)$ does not lie on a $k$-plane
generated by $k+1$ points in $X$, then $\rg \varphi \geq d-k+1$. 
\end{prop}
\begin{proof}
We know by lemma \ref{secante} that $\varphi$ satisfies the condition
\begin{equation}\label{uno}
\varphi\left(L_1^{n_1}\dots L_r^{n_r}.f\right)=0 
\hskip 20pt \forall f\in S_{d-k-1}
\end{equation} 
where $L_1,\dots,L_r$ are distinct linear forms, and $\sum_{i=1}^r n_i =k+1$. 
As $[\varphi]$ does not lie
on a secant $k$-plane  to $X$,
by using lemma \ref{igual} we can assume that $n_1\geq 2$. 
Let $[\beta_i] \in \zP^1(K)$ be the zero of $L_i$ for $i=1,\dots,r$.

Assume that $\varphi$ is a linear combination of $d-k$ elements in $X$,
\begin{equation}\label{dos}
\varphi=\sum _{j=1} ^{d-k} \lambda_j ev(\alpha_j).
\end{equation}
Since $ r \leq k$, and $k\leq\left[\frac {d-1} 2 \right]$,
then $r<d-k$. Therefore at least one $[\alpha_j]$
differs from all $[\beta_i]$.
We can assume  that
$[\beta_i] = [\alpha_i]$ for all $i \leq s$,
and $[\beta_i] \neq [\alpha_j]$ for all $i,j>s$, 
where $s\leq r$.
Consider for every $j$, $1 \leq j \leq d-k$,
a linear form $M_j$ having $\alpha_j$ as its zero, 
and let $g$ be the following polynomial of degree $d$ 
$$
g=\left(\prod _{j=1}^s  M_j \right) \left(\prod _{j=s+2}^{d-k}M_j
\right) \left( \prod_{i=1}^rL_i^{n_i}\right)
$$
We have $\varphi(g)=0$ by equation (\ref{uno});
and $\varphi(g)=\lambda_{s+1}.c$, where $c \neq 0$, by equation (\ref{dos}).
Then $\lambda_{s+1}=0$.
Similarly, we get $\lambda_j=0$ for all $j\geq s+1$.

Therefore $\varphi =\sum _{j=1} ^{s} \lambda_j ev(\beta_j)$, 
and so $\varphi\left(L_1\dots L_s.f\right)=0$ for all $f\in S_{d-s}$.
But this means that $[\varphi]$ lies on an $(s-1)$-plane  generated by 
$s$ points in $X$, which contradicts the fact that $[\varphi]$
does not lie on a secant $k$-plane.
\end{proof}

\bigskip

\begin{lema}\label{vuelta}
We consider for each integer $k\leq \left[\frac {d-1} 2\right]$ 
the following vector subspace of $S_{d-k}$
$$
A_k=\{f\in S_{d-k}: \varphi(fg)=0\;\forall g\in S_k\}.
$$ 
Then
   \begin{enumerate}
   \item
   $\codim (A_k)=r \leq k$ if and only if  $[\varphi]\in \se^{r-1}(X)$.
   \item If $ \codim(A_k)=k+1$ and $A_k\subset\Delta_{d-k}$, then 
   $[\varphi]\in \se^k(X)\setminus \se^{k-1}(X)$ 
   and does not lie on a secant $k$-plane. 
   \end{enumerate}
\end{lema}

\begin{proof}
The first statement follows from lemma \ref{secante}.

Let $A=A_k$.
For the second statement, we use Bertini's Theorem to 
show that $A$ has base points,
since all polynomials  in $A$ have singular points.
Moreover, one of these base points is a multiple root
for every polynomial in $A$.

Let $L_1,L_2,\dots,L_r$ be linear forms having as zeros the base points of $A$.
We have $A=L_1^{n_1}L_2^{n_2}\dots
L_r^{n_r}Z$, where $\sum_{i=1}^r n_i=m\leq k+1$
and $Z\subseteq S_{d-k-m}$ is 
a linear subspace of dimension $d-2k$ having no base points.
We claim that $m=k+1$. If $m < k+1$, then $Z$ is a proper subspace
of $S_{d-k-m}$ having no base points. Therefore, 
$\dim S_1.Z\geq \dim Z+2$  (\cite{Ha}, lemma 9.8 pag. 103),
and the codimension of  $S_1.Z$ in  $S_{d-k-m+2}$ is less than or equal 
to $k-m<k-m+1=\codim(Z)$.
If $S_1.Z$ is a proper subspace,
we can use the same result since  $S_1.Z$ will have no base points.
Repeating this procedure at most $k-m+1$ times,
we get a space of codimension 0,
therefore $S_{k-m+1}.Z=S_{d-2m+1}$.
We conclude that $S_{k-m+1}.A=L_1^{n_1}L_2^{n_2}\dots
L_r^{n_r} S_{k-m+1}Z=L_1^{n_1}L_2^{n_2}\dots
L_r^{n_r}S_{d-2m+1}$.
Then 
$$
0=\varphi(S_{k-m+1}.A.S_{m-1})=\varphi(L_1^{n_1}\dots L_r^{n_r}
S_{d-2m+1}S_{m-1})=\varphi(L_1^{n_1}\dots L_r^{n_r}
S_{d-m}).
$$
This would mean that $[\varphi] \in \se^{m-1}(X)$, 
which contradicts the fact that  $\codim (A)=k+1$.

Therefore $k+1=m$, $A=L_1^{n_1}\dots L_r^{n_r}S_{d-k-1}$,
and $[\varphi] \in \se^k (X)$. 
Since $A \subset \Delta_{d-k}$, we know by lemma \ref{igual} that $[\varphi]$ 
does not lie on a secant $(d-k-1)$-plane. Then $[\varphi]$ 
does not lie on a secant $k$-plane, because $k\leq d-k-1$. 
\end{proof}

\medskip

The next proposition is  reciprocal to proposition \ref{ida}. 
\medskip
\begin{prop}\label{vuelta2}
Let $k\leq \left[\frac {d-1} 2\right]$, and let 
$\varphi$ be a form  having rank greater or equal to $d-k+1$,
$[\varphi] \not \in \se^{k-1}(X)$.
Then $[\varphi] \in \se^k (X)\setminus \se^{k-1} (X)$ 
and $[\varphi]$ does not lie on a secant $k$-plane.
\end{prop}
\begin{proof}
It follows directly from lemmas \ref{igual} and \ref{vuelta}.
\end{proof}

\medskip

Now we can prove  theorem \ref{elteorema}.

\medskip

\begin{proof}[Proof of theorem \ref{elteorema}]
We rewrite the  statement in theorem 2 using secant varieties:
$$\se^k(X) \setminus \se^{k-1}(X)
=\{[\varphi] : \rg \varphi = k+1\} \cup 
\{[\varphi]  : \rg \varphi = d-k+1\}
$$

Note that this  implies
$$\se^k(X) 
=\{[\varphi] : \rg \varphi\leq k+1\} \cup \{[\varphi]: \rg \varphi \geq d-k+1\}
$$
for each $k$.

We will use induction on $k$.

For  $k=0$, $\se^0(X)=X=\{[\varphi]: \rg \varphi=1\} \cup \emptyset$.

For $k>0$, let  $[\varphi] \in \se^k(X) \setminus \se^{k-1}(X)$.
If $[\varphi]$ lies on a secant $k$-plane, its rank is by definition $k+1$.
Otherwise, we have $\rg \varphi \geq d-k+1$
by proposition \ref{ida}.
If $\rg \varphi > d-k+1$, $[\varphi]$
lies on $\se^{k-1}(X)$ by induction, which contradicts the fact
that $[\varphi]\in \se^k(X) \setminus \se^{k-1}(X)$.
Therefore $\rg (\varphi)=d-k+1$. 

We have proved that 
$$\se^k(X) \setminus \se^{k-1}(X)
\subseteq \{[\varphi] : \rg \varphi=k+1\} \cup \{[\varphi]: \rg \varphi=d-k+1\}.
$$

If $\rg \varphi =k+1$, then, by definition, $[\varphi] \in \se^k(X)$; and
 $[\varphi]\notin\se^{k-1}(X)$ by induction.
Finally, if $\rg(\varphi)=d-k+1$ then $[\varphi] \not \in
\se^{k-1}(X)$ by induction and  $[\varphi] \in \se^k(X)$ by proposition \ref{vuelta2}.

In the case $d=2n+1$, this procedure ends when $k=n$,
and we have
$$\se^n(X) \setminus \se^{n-1}(X)
=\{[\varphi] : \rg \varphi=n+1\} \cup \{[\varphi]: \rg \varphi=n+2\}.
$$

In the case $d=2n$ and $k=n$ we cannot use proposition \ref{vuelta2}.
However, we know that $\se^{n-1}(X)
= \{[\varphi] : \rg \varphi \leq n\} \cup \{[\varphi]: \rg \varphi \geq n+2\}$,
and then every form not lying on $\se^{n-1}(X)$
has rank $n+1$. Therefore, since $\se^{n}(X)= \zP(S_d^\ast)$,
we have $\se^{n}(X)\setminus \se^{n-1}(X)
=\{[\varphi]: \rg \varphi =n+1\}$.

\end{proof}

\end{section}
\begin{section}{Computing the rank}
Next we show how to determine the rank of a binary form $q\in S_d$.
Firstly we  find the integer $k$ 
such that 
$q \in \overline{S_{d,k+1}} \setminus \overline{S_{d,d-k+1}}$.
Indeed, if $q=Z_0x^d+\binom{d}{1}Z_1x^{d-1}y+\dots+\binom{d}{d-1}Z_{d-1}xy^{d-1}+Z_dy^d$
then $k+1$ is the rank of one of these two matrices
$$
\left [
\begin{array}{cccc}
Z_0 & Z_1 & \dots & Z_{n}\\
Z_1 & Z_2 & \dots & Z_{n+1}\\
\vdots & \vdots &  & \vdots\\
Z_{n} & Z_{n+1} & \dots & Z_d\\
\end{array}
\right] \hskip 20pt \text{ or } \hskip 20pt 
\left [
\begin{array}{cccc}
Z_0 & Z_1 & \dots & Z_n\\
Z_1 & Z_2 & \dots & Z_{n+1}\\
\vdots & \vdots &  & \vdots\\
Z_{n+1} & Z_{n+2} & \dots & Z_d\\
\end{array}
\right]
$$
depending on wether $d=2n$ or $d=2n+1$ (by proposition \ref{harris}).

Secondly, we decide
if $\rg q=k+1$ or $\rg q=d-k+1$.
If $d=2n$, and $k=n$ then $k+1=d-k+1=\rg q$.
Otherwise, let
$$
M=\left [
\begin{array}{cccc}
Z_0 & Z_1 & \dots & Z_{k+1}\\
Z_1 & Z_2 & \dots & Z_{k+2}\\
\vdots & \vdots &  & \vdots\\
Z_{d-k-1} & Z_{d-k} & \dots & Z_d\\
\end{array}
\right]
$$
be the matrix of the bilinear form $B:S_{d-k+1} \times S_{k+1}
\rightarrow K$ defined in the proof of proposition \ref{secante}.
If $\varphi$  denotes the linear functional dual to $q$,
then $B$ is the bilinear form associated to $\varphi$.
The rank of this matrix is also $k+1$ (by proposition \ref{harris}),
so there is a polynomial $f \in S_{k+1}$, unique up to scalars,
such that $B(S_{d-k+1},f)=0$.

In order to find $f$, we solve the system
$$
M\cdot\begin{pmatrix} f_0 \\ f_1 \\ \vdots \\ f_k \\ f_{k+1} \end{pmatrix} =0
$$
and set  $f=f_0x^{k+1}+f_1x^ky+\dots+f_kxy^k+f_{k+1}y^{k+1}$.

If $f$ has no multiple roots, then $\rg q=k+1$, otherwise
$\rg q=d-k+1$. Indeed, if $f$ has $k+1$ distinct
roots $[t_1,u_1],\dots,[t_{k+1},u_{k+1}] \in \zP^1(K)$, 
then $q$ is a linear combination of the $d$th powers
of the linear forms $L_1=t_1x +u_1y,\dots ,L_{k+1}=t_{k+1}x + u_{k+1}y$.

\end{section}


\begin{thebibliography}{99}
\bibitem[Ha]{Ha}Harris J.,
{\it Algebraic Geometry},
Graduate Texts in Mathematics,
Springer Verlag, Corrected third prining 1995.
\end{thebibliography}
\end{document}